\documentclass[a4paper]{article}
\usepackage{amsmath}
\usepackage{breqn}
\usepackage{amsfonts}
\usepackage{amsthm}
\usepackage{amssymb}
\usepackage{hyperref}
\usepackage{graphicx} 
\usepackage{booktabs}
\usepackage[T1]{fontenc}
\usepackage[utf8]{inputenc}
\usepackage{authblk}
\newtheorem{theorem}{Theorem}[section]
\newtheorem{corollary}[theorem]{Corollary}
\newtheorem{lemma}[theorem]{Lemma}

\newtheorem{definition}[theorem]{Definition}

\title{ Stochastic ordering results in parallel and series systems with  Gumble distributed random variables}

\author[1]{Surojit Biswas\thanks{$^\ast$ Corresponding author. e-Mail: sb38@iitbbs.ac.in}}
\author[2]{Nitin Gupta\thanks{$^\ast$ Co-author. e-Mail: nitin.gupta@maths.iitkgp.ernet.in}}
\affil[1,2]{Department of Mathematics, Indian Institute of Technology Kharagpur}

\begin{document}

\maketitle

\begin{abstract}
\textbf{The stochastic comparisons of parallel  and series system are worthy of study. In this paper, we present some stochastic comparisons of parallel and series systems having  independent components from Gumble distribution with two parameters (one location and one shape).  Here, we first put a  condition for the likelihood ratio ordering of the parallel systems and second we use the concept of vector majorization technique to compare the systems by the reversed hazard rate ordering, the hazard rate ordering, the dispersive ordering, and the less uncertainty ordering with respect to the location parameter.}
\end{abstract}
\section{Introduction}
Discussion about order statistics is quite important for different distributions. A wide range of application of order statistics can be found in many different areas like   statistics, applied probability, reliability theory, actuarial science, auction theory, hydrology, life testing \textit{etc}. \\

  Let $ X_{1},.....,X_{n} $ be a set of independently distributed random variables and $ X_{n:n}=max\left\lbrace X_{1},.....,X_{n}\right\rbrace, X_{1:n}=min\left\lbrace X_{1},.....,X_{n}\right\rbrace$. 
 $ X_{1:n} $  represents a series system and $ X_{n:n} $  represents a parallel system. $ X_{1:n} $ and $ X_{n:n} $ are also known as $ 1$st and $ n $th order statistics respectively. It is well known that in a parallel system at least one of it's $ n $ components needs to work and in a series system all the components of the system need to work at a time. That's why the parallel and the series system are called as $ 1$-out-of-$n $ and $ n$-out-of-$n $ systems respectively. Thus, the study of stochastic comparisons of $ m $-out-of-$ n $ systems is similar to the study of the stochastic comparisons of the order statistics.\\
 
 We use extreme value theory to study  the stochastic behavior of extreme values of a process.  The Gumble distribution is one of the extreme value distribution which is also known as an extreme value type-I or smallest extreme value (Type I) distribution due to Emil Julius Gumbel(1954)\cite{p9}. We can observe that the maximum of a random sample after proper randomization can only converge to  one of the three possible
 distributions, the Gumbel distribution, the Fr\'{e}chet distribution or the Weibull distribution (see \cite{p8})\\

  Various results on stochastic comparisons of the Weibull distribution and Fr\'{e}chet distribution  are already presented by  Torrado et. al.(see \cite{p4}), Fang et. al.(see \cite{p5}) and Gupta et. al. (see \cite{p1}). In this paper, we present several results for the parallel and series systems having independent Gumble distributed components. The cumulative distribution function (cdf) of the Gumble distribution ($ ' $Gum$ ' $ in short) is given by
\begin{equation}
F(x)= e^{-e^{-\frac{(x-\mu)}{\sigma}}}
\hspace{.4cm} x,\mu\in \mathbb{R}, \sigma>0.
\end{equation}
 Here $ \mu $ is the location parameter and $ \sigma $ is the scale parameter. A random variable $ X $ is said to follow Gumble distribution, denoted as  $ X\sim Gum(\mu,\sigma) $ if $ X $ has the cumulative distribution function given in (1). It is very much useful for analyzing  several extreme natural events like earthquakes, floods, rainfall, wind speeds, snowfall \textit{etc}. For applications
 of the distribution distribution one may refer to (\cite{p9},\cite{p10}).\\

 Many researchers studied  order statistics in terms of  stochastic comparisons. A vast literature on stochastic comparisons for the lifetimes of the series and parallel systems are already available where the component variables follow  Fr\'{e}chet\cite{p1}, Generalized Exponential\cite{p6}, Exponentiated Gamma\cite{p12}, Exponentiated Scale model\cite{p13}, Exponential Weibull\cite{p14} distributions \textit{etc}. The comparisons are made of with respect to the usual stochastic order, hazard rate order, reversed hazard rate order, likelihood ratio order, dispersive order, \textit{etc}. For further details on stochastic
 comparisons, one may refer to (\cite{p4},\cite{p5},\cite{p6},\cite{p7},\cite{p8}).
  See, Shaked M et. al.\cite{b2} for more comprehensive discussion on order statistics.\\
  
  The aim of this paper is to present the likelihood ratio, reversed hazard rate, hazard rate, dispersive, and the less uncertainty  ordering for parallel, series systems having Gumble distributed components. 
 Now, let $ X_{1},.....,X_{n} $ be continuous independent random variables such that, $ X_{i}\sim Gum(\mu_{i}, \sigma), \hspace{.3cm} i=1,2,....,n$. Furthermore, let $ Y_{1},.....,Y_{n} $ be another set of continuous independent random variables such that $ Y_{i}\sim Gum( \mu_{i}^{*}, \sigma),  i=1,2,....,n$. When  $ \mu_{i}\geq \mu^{*}_{i} ,$ we obtain the likelihood ratio ordering for the parallel systems, and when  $ (\mu_{1},\cdotp\cdotp\cdotp\cdotp\cdotp\mu_{n})\succeq^{m}(\mu^{*}_{1},\cdotp\cdotp\cdotp\cdotp\cdotp\mu^{*}_{n}) ,$ we proceed as follows.
 First, we discuss the  reversed hazard rate ordering for the parallel systems.
 Second, we discuss the hazard rate ordering for  the series systems. Finally, we discuss the dispersive and less uncertainty ordering for the series systems. \\
   
   The paper has been organized in the following manner:\\
   In Section-2 we briefly present some basic useful definitions, lemmas, and theorems which we have used throughout this paper. In the Section-3,  we deal with the concept of vector majorization technique and achieve several ordering results. 
   \section{Preliminaries}
This section consists of some important definitions, theorems, and lemmas  that are most pertinent to developments in Section-3..  We use the notation  $\mathbb{R}$=$(-\infty,+\infty),$  $ \mathbb{R}^{+}$=$[0,+\infty),$  $\textquoteleft\log$' for usual logarithm base $ e,$ and we use $ \textquoteleft $increasing' and $ \textquoteleft$decreasing' for $ \textquoteleft $non-decreasing' and $ \textquoteleft $non-increasing' respectively throughout this paper.\par
   
   Let $ X, Y $ be two continuous  random variables having the cumulative distribution 
   functions $ F_{X}(\cdot)$ and $ F_{Y}(\cdot)  ,$ the survival functions $ \bar{F}_{X}(\cdot)$ and $ \bar{F}_{Y}(\cdot),$ density functions $ f_{X}(\cdot) $ and $ f_{Y}(\cdot) ,$ hazard rate functions $ r_{X}(\cdot)=\dfrac{f_{X}(\cdot)}{\bar{F}_{X}(\cdot)} $ and $ r_{Y}(\cdot)=\dfrac{f_{Y}(\cdot)}{\bar{F}_{Y}(\cdot)} ,$ and reversed hazard rate functions $ \tilde{r}_{X}(\cdot)=\dfrac{f_{X}(\cdot)}{F_{X}(\cdot)} $ and  $\tilde{r}_{Y}(\cdot)=\dfrac{f_{Y}(\cdot)}{F_{Y}(\cdot)} .$\\

   Let $ X $ be a continuous random variable having $ F_{X}(\cdot)$ and  $ f_{X}(\cdot)$ as cumulative distribution and density functions respectively. The  measure of entropy of $ X $ due to Shannon(1948)\cite{p14} is defined by
   $$ \mathcal{H}(f_{X})=-E\left[ \log f_{X}(x)   \right]=-\int_{0}^{\infty}f_{X}(x)\log f_{X}(x), $$ which is commonly known as Shannon information measure. This is a measure of the uncertainty of the lifetime of a system.  A system having  low uncertainty  is more reliable than a system with great uncertainty. Ebrahimi and Pellerey (1995)\cite{p2} defined the uncertainty of residual lifetime distributions, $ \mathcal{H}(f_{X},t) $, of a component by
$$ \mathcal{H}(f_{X},t)=\mathcal{H}(X-1|X\geq t)=-\int_{t}^{\infty}\dfrac{f_{X}(x)}{\bar{F}(t)} \log \dfrac{f_{X}(x)}{\bar{F}(t)}  dx  $$
$$\hspace{5.5cm}=\log \bar{F}(t)-\frac{1}{\bar{F}(t)}\int_{t}^{\infty} f_{X}(x) \log f_{X}(x)  dx$$
$$\hspace{4.5cm}=1-\frac{1}{\bar{F}(t)}\int_{t}^{\infty} f_{X}(x) \log r_{X}(x)  dx,$$ where $ r_{X}(\cdot), \bar{F}(\cdot) $ are hazard and survival functions respectively. 
After the component has survived upto time $ t $, $ \mathcal{H}(f_{X},t) $ measures the expected uncertainty contained in the
conditional density of $ X - t $ given $ X> t $ about the predictability of the remaining
lifetime of the component. Now for $ t=0, $ we can see that $$ \mathcal{H}(f_{X},0)=-\int_{0}^{\infty}f_{X}(x)\log f_{X}(x), $$ which is nothing but the Shannon's  measure of entropy of $ X $.\\

Next, we present some useful definitions for different stochastic ordering.

\begin{definition}
 (\textit{Stochastic Order}) \\
 Let $ X, Y $ be two random variables. We say  that $X$ is smaller than $Y$ in
 \begin{enumerate}
 	 	\item the \textit{likelihood ratio order} (denoted by, $X \leq_{lr} Y$) if $ \dfrac{f_{Y}(x)}{f_{X}(x)}$ is increasing in $ x $.
 	\item the \textit{reversed hazard rate order} (denoted by, $X \leq_{rh} Y$) iff $\tilde{r}_{X}(x) \leq \tilde{r}_{Y}(x),$ $ x \in \mathbb{R}.$ Equivalently  $ \dfrac{F_{Y}(x)}{F_{X}(x)}$ is increasing in $x$.
 
 	\item  the \textit{hazard rate order} (denoted by, $X \leq_{hr} Y$)
 	iff $r_{Y}(x) \leq r_{X}(x),\hspace{.2cm} x \in \mathbb{R}.$ Equivalently if $ \dfrac{\bar{F}_{Y}(x)}{\bar{F}_{X}(x)}$ is increasing in $x$.

 	\item  the \textit{dispersive} order (denoted by, $X \leq_{disp} Y$ ) if for $ 0\leq \alpha<\beta\leq1 $ we have
 	$$ F^{-1}_{X}(\beta)-F^{-1}_{X}(\alpha)\leq F^{-1}_{Y}(\beta)-F^{-1}_{Y}(\alpha).            $$ Equivalently $X \leq_{disp} Y$ iff $f_{Y}(F^{-1}_{Y}(p))\leq f_{X} (F^{-1}_{X}(p))$ for all $ p\in (0,1). $
 	\item the less uncertainty order (denoted by $X \leq_{LU} Y$ ) if $ \mathcal{H}(f_{X};t)\leq \mathcal{H}(f_{Y};t) $ for any $ t>0. $
 	
 \end{enumerate}
The following relation is well-known, that is \\

$$   X \leq_{lr} Y \Longrightarrow X \leq_{hr} Y(X \leq_{rh} Y)\Longrightarrow X \leq_{st} Y.$$ 
	
\end{definition}
\begin{definition}
\textit{(Majorization)}\\
Let $ \textbf{u}=(u_{1},...,u_{n})$ and $ \textbf{v}=(v_{1},...,v_{n})$   be two real vectors from $ \mathbb{R}^{n} $ with the order components, $ u_{(n)}\leq.....\leq u_{(1)} $ and $v_{(n)}\leq.....\leq v_{(1)}, $ respectively. Then we say $ \textbf{u}$ majorizes $\textbf{v} ,$ denoted by
$ \textbf{u}\succeq^{m}\textbf{v}$ if $$ \sum_{i=1}^{k}u_{i} \leq \sum_{i=1}^{k}v_{i}$$ $ k=1,2,...,n-1,$ and \hspace{.04cm} $\sum_{i=1}^{n}u_{i}= \sum_{i=1}^{n}v_{i}.$ 
\end{definition}
\begin{definition}
	Let $ \textbf{u}=(u_{1},.....,u_{n})$ and $ \textbf{v}=(v_{1},.....,v_{n})$ be two vectors from $ \mathbb{R}^{n} $. A real valued function $ \sigma(\textbf{u}):\mathbb{R}^{n}\rightarrow \mathbb{R} $ is said to be \textit{Schur-concave} and \textit{Schur-convex}  if  for all $ \textbf{u}\succeq^{m} \textbf{v} $ we have $ \sigma(\textbf{u})\leq \sigma(\textbf{v})$ and $\sigma(\textbf{u})\geq \sigma(\textbf{v}),$ respectively.
\end{definition}\par
The theorem, stated bellow is very useful for our results.
\begin{theorem} \textit{(Marshall et al., p.84, \cite{b1})}:
	Let $ \mathcal{I}\subset\mathbb{R} $ be an
	open interval and let $ \sigma: \mathcal{I}^{n} \rightarrow \mathbb{R}$ be continuously differentiable function. The necessary and sufficient conditions for $ \sigma $ to be Schur-convex(Schur-concave) on $ \mathcal{I}^{n}$ are
	$ \sigma $ is symmetric on $ \mathcal{I}^{n}$ and, for all $ i\ne j $
	\begin{center}
	$ (z_{i}-z_{j})\left( \dfrac{\partial \sigma}{\partial z_{i}}(\textbf{z})-\dfrac{\partial \sigma}{\partial z_{j}}(\textbf{z})\right)\geq 0 \hspace{.1cm} (\leq 0) $
	\end{center}	for all $ \textbf{z}\in \mathcal{I}^{n} $. Where, $  \dfrac{\partial \sigma}{\partial z_{i}}$ is partial derivative of $ \sigma $ with respect to the $ i^{th} $ component of $ \textbf{z} $.
	
\end{theorem}

\begin{lemma}
	Let the function $ \phi(t):\mathbb{R^{+}}\rightarrow
	 \mathbb{R^{+}}$ be defined as
	 \begin{center}
	  $ \phi(t)=\dfrac{t}{e^{t}-1}. $
	 \end{center}
 Then, 
 \begin{enumerate}
 	\item $ \phi(t) $ is a convex function in $ \mathbb{R^{+}}; $
 	\item  $ \phi(t)  $ is decreasing with respect to $ t .$
 \end{enumerate}
\begin{proof}
	\begin{enumerate}
		\item The proof of this lemma is already avaiable in \cite{p1}, and therefore skipped for shake of brevity.
		\item Taking derivative of $ \phi(t) $ with respect to $ t $ we have
		$$ (e^{t}-1)^{2}\dfrac{d\phi}{dt}=e^{t}-1-te^{t}.         $$ Now, let $ g(t)=e^{t}-1-te^{t}, $ then $ \dfrac{dg}{dt}= e^{t}-e^{t}-te^{t}.$ Which implies that $\dfrac{dg}{dt}=-te^{t}\leq0$ for any $ t\in \mathbb{R^{+}}. $ So, $ g(t)\leq g(0)=0.$ Therefore, we have $ \dfrac{d\phi}{dt}\leq 0,$ this means that $ \phi(t) $ is decreasing in $ t. $
	\end{enumerate}
\end{proof}

\end{lemma}
\begin{lemma}
	Let $ \mathcal{I}\subset \mathbb{R} $ be an open interval such that if a function $ \gamma(x):\mathcal{I}\rightarrow\mathbb{R}
	$ is convex then the function $h(\boldsymbol{x}) $ defined as
   $$h(\boldsymbol{x})=\sum_{i=1}^{n}\gamma(x_{i}) $$
is a Schur-convex function in $ \mathcal{I}^{n}, $ where $ \boldsymbol{x}=(x_{1},...,x_{n}). $ Consequently for $ \boldsymbol{x}=(x_{1},...,x_{n}), $ $ \boldsymbol{y}=(y_{1},...,y_{n})\in \mathcal{I}^{n} $ if $ (x_{1},...,x_{n})\succeq^{m}(y_{1},...,y_{n})$ then $ h(\boldsymbol{x})\geq h(\boldsymbol{y}). $
	\begin{proof}
	Proof of this lemma can be found in Marshall et al.\cite{b1}
\end{proof}
\end{lemma}

\begin{theorem}
	Let $ X,Y $ be two independent random variables having density functions $ f_{X}, f_{Y} $, cumulative distribution functions $ F_{X}, F_{Y} $, survival functions $ \bar{F}_{X}, \bar{F}_{Y} ,$ respectively. Let $ T $ be a random variable with density function $ h $ and distribution function $ H $. $ T $ is independent of $ X $ $ Y $. Then if
	\begin{enumerate}
		\item 	$ X\geq_{rh} Y  $ and either $ X $ or $ Y $ is IRHR then $ X^{T}\geq_{rh} Y^{T}. $ 
		\item 	$ X\leq_{hr} Y  $ and either $ X $ or $ Y $ is DHR then $ X^{T}\leq_{hr} Y^{T}. $
	\end{enumerate}
\begin{proof}
	The proof is similar to the proof in \cite{p3}.
\end{proof}
\end{theorem}
Sometimes it may not be possible to find dispersive ordering directly from the definition. In some particular cases, we can use the following theorem to have dispersive ordering.
\begin{theorem}
	Let $ X $ and $ Y $ be two  random variables. If $ X\leq_{hr}Y $ and $ X $ or $ Y $ is DHR, then  $ X\leq_{disp}Y. $
	\begin{proof}
		For the proof see Shaked M et. al.\cite{b2}
	\end{proof}
\end{theorem}
The result of the Theorem-2.3 from Ebrahimi et.al.\cite{p2} can be strengthened for the hazard rate ordering as follows 
\begin{theorem}
	Let $ X $ and $ Y $ be two random variables. If $ X\leq_{hr}Y $ and $ X $ or $ Y $ is DHR, then  $ X\leq_{LU}Y. $
	\begin{proof}
		Proof can be found in \cite{p2}.
	\end{proof}
	
\end{theorem}

\section{ Results}
  
 In this section, we work on several results by comparing the lifetimes of parallel and series systems having independent Gumble distributed components. These results are presented for the location parameter $ \mu $ by using vector majorization technique.\\
  
 In the following theorem, we put a sufficient condition on the location parameter $ \mu $ for the likelihood ratio ordering of the parallel systems where the components of the systems follow the Gumble distribution. 
 \begin{theorem}
 	Let $ X_{1},....,X_{n}$ $ (Y_{1},....,Y_{n})$  be the set of continuous independent random variables having $ X_{i}\sim Gum(\mu_{i}, \sigma)(Y_{i}\sim Gum(\mu^{*}_{i}, \sigma)),$ $ i=1,2,....,n$. Then, if  $ \mu_{i}\geq \mu^{*}_{i} $, we have
 	\begin{center}
 		$ X_{n:n}\geq_{lr} Y_{n:n}$.
 	\end{center}
 \begin{proof}
 	The cumulative distribution function of $ X_{n:n} $  is given by 
 	$$F_{X_{n:n}}(x)=\prod_{i=1}^{n} e^{-e^{-\frac{(x-\mu_{i})}{\sigma}}}\hspace{.2cm}\hspace{.4cm} x,\mu_{i}\in \mathbb{R}, \sigma>0, i=1,...,n. 
 	 $$ 
 	Therefore, the density function of $ X_{n:n}$ is
 		$$ f_{X_{n:n}}=\dfrac{d}{dx}\prod_{i=1}^{n} e^{-e^{-\frac{(x-\mu_{i})}{\sigma}}}$$ $$\hspace{3.66cm}= \dfrac{F_{X_{n:n}}(x)}{\sigma}\sum_{i=1}^{n}e^{-\frac{(x-\mu_{i})}{\sigma}}, \hspace{.3cm}i=1,...,n.$$ 
 		Similarly, the density function of $ Y_{n:n}$ is
 		$$f_{Y_{n:n}}= \dfrac{F_{Y_{n:n}}(x)}{\sigma}\sum_{i=1}^{n}e^{-\frac{(x-\mu^{*}_{i})}{\sigma}}, \hspace{.3cm}i=1,...,n.$$

Now, Let us consider $$ lr(x)=\frac{f_{X_{n:n}}(x)}{f_{Y_{n:n}}(x)}= \dfrac{F_{X_{n:n}}(x)}{F_{Y_{n:n}}(x)} \left[ \dfrac{\sum_{i=1}^{n}e^{-\frac{(x-\mu_{i})}{\sigma}}}{\sum_{i=1}^{n} e^{-\frac{(x-\mu^{*}_{i})}{\sigma}}}\right] $$ 
$$\hspace{3cm}=\left[ \dfrac{\sum_{i=1}^{n}e^{\frac{\mu_{i}}{\sigma}}}{\sum_{i=1}^{n} e^{\frac{\mu^{*}_{i}}{\sigma}}}\right]  \prod_{i=1}^{n}e^{-\left[ e^{-\frac{(x-\mu_{i})}{\sigma}}-e^{-\frac{(x-\mu^{*}_{i})}{\sigma}}\right] },\hspace{.3cm}i=1,..,n. $$ Taking derivative with respect to $ x $ we get 
 $$ \dfrac{d}{dx}lr(x)=\frac{1}{\sigma} \dfrac{F_{X_{n:n}}(x)}{F_{Y_{n:n}}(x)} \left[ \dfrac{\sum_{i=1}^{n}e^{\frac{\mu_{i}}{\sigma}}}{\sum_{i=1}^{n} e^{\frac{\mu^{*}_{i}}{\sigma}}}\right]  \sum_{i=1}^{n}\left[ e^{-\frac{(x-\mu_{i})}{\sigma}}-e^{-\frac{(x-\mu^{*}_{i})}{\sigma}}\right],\hspace{.3cm}i=1,..,n. $$ Now, by our assumption $ \mu_{i}\geq \mu^{*}_{i} $ and since $e^{-\frac{(x-\mu)}{\sigma}}  $ is increasing in $ \mu ,$ so, we have $e^{-\frac{(x-\mu_{i})}{\sigma}}\geq e^{-\frac{(x-\mu^{*}_{i})}{\sigma}}, i=1,...,n. $ Using this observation and since all the  other terms in the right hand side of the above equation is always garter than 0. So, we finally conclude  that $ \dfrac{d}{dx}lr(x)\geq0, $  which implies that $ lr(x) $ is increasing in $ x $ and this proves our theorem.

 \end{proof}
 
 \end{theorem}

  The following theorem deals with the reversed hazard rate ordering of a parallel system. The components of the system are having fixed scale parameter $ \sigma $ and a varying location parameter $ \mu $ of the Gumble distribution.
\begin{theorem}
Let $ X_{1},....,X_{n}$ $ (Y_{1},....,Y_{n})$  be the set of continuous  independent random variables having $ X_{i}\sim Gum(\mu_{i}, \sigma)(Y_{i}\sim Gum(\mu^{*}_{i}, \sigma)),$ $ i=1,2,....,n$. Then, if  $ (\mu_{1},\cdotp\cdotp\cdotp\cdotp\cdotp\mu_{n})\succeq^{m}(\mu^{*}_{1},\cdotp\cdotp\cdotp\cdotp\cdotp\mu^{*}_{n}) $, we have
		 \begin{center}
		$ X_{n:n}\geq_{rh} Y_{n:n}$.
	\end{center}
  	
\begin{proof}
	
	 The  cumulative distribution function of $ X_{n:n} $ is
	$$F_{X_{n:n}}(x)=\prod_{i=1}^{n} e^{-e^{-\frac{(x-\mu_{i})}{\sigma}}}
	\hspace{.4cm} x,\mu_{i}\in \mathbb{R}, \sigma>0, \hspace{.1cm} i=1,...,n. $$

Now, we know that, in a parallel system, the sum of the reversed hazard rate of the lifetime of each
	component is equal with the reversed hazard rate of the lifetime of the system. 
	Therefore, the reversed hazard rate function of $ X_{n:n}$ is 
	
	$$ \tilde{r}_{X_{n:n}}=\frac{1}{\sigma}\sum_{i=1}^{n}e^{-\frac{(x-\mu_{i})}{\sigma}},\hspace{.2cm} i=1,...,n.$$ It is clear that $ e^{-\frac{(x-\mu_{i})}{\sigma}}$ is convex in $ \mu_{i}$ for $i=1,...,n. $ Therefore, from Lemma-2.6 we obtain $\sum_{i=1}^{n}e^{-\frac{(x-\mu_{i})}{\sigma}}$ is Schur-convex with respect to $ \mu_{i}, i=1,...,n, $ which means that $ \tilde{r}_{X_{n:n}}\geq \tilde{r}_{Y_{n:n}}$ i.e. $X_{n:n}\geq_{rh} Y_{n:n}, $ as desired.
	
\end{proof}
\end{theorem}
The above theorem leads us to the following corollary
\begin{corollary}
	Let $ X_{1},....,X_{n}$ $ (Y_{1},....,Y_{n})$  be the set of continuous  independent random variables having $ X_{i}\sim Gum(\mu_{i}, \sigma)(Y_{i}\sim Gum(\mu_{i}, \sigma)),$ $ i=1,2,....,n$. Then, if  $ (\mu_{1},\cdotp\cdotp\cdotp\cdotp\cdotp\mu_{n})\succeq^{m}(\mu^{*}_{1},\cdotp\cdotp\cdotp\cdotp\cdotp\mu^{*}_{n}) ,$ and let $ T $ be any random variable which is independent of $ X $ and $ Y $. Then we have
	\begin{center}
		$ X^{T}_{n:n}\geq_{rh} Y^{T}_{n:n}$.
	\end{center}

\begin{proof}
	To prove $ X^{T}_{n:n}\geq_{rh} Y^{T}_{n:n}$. By Theorem-2.7(1), it is sufficient if we can show $ X_{n:n}\geq_{rh} Y_{n:n}$ and $ X_{n:n} $ is IRHR. Now, we know the reversed hazard rate of $ X_{n:n} $ with respect to $ \mu_{i}, i=1,...,n, $ is 
		$$ \tilde{r}_{X_{n:n}}=\frac{1}{\sigma}\sum_{i=1}^{n}e^{-\frac{(x-\mu_{i})}{\sigma}},\hspace{.2cm} i=1,...,n.$$ 
		Therefore, it is easy to cheack that $\sum_{i=1}^{n}e^{-\frac{(x-\mu_{i})}{\sigma}} $ is increasing in $ \mu_{i}, i=1,...,n. $ So, $ \tilde{r}_{X_{n:n}} $ is IRHR. Using this observation and the result in Theorem-3.2, we get $ X^{T}_{n:n}\geq_{rh} Y^{T}_{n:n},$ which completes the proof.
\end{proof}

\end{corollary}

Next theorem discusses the comparison of the series system with respect to hazard rate ordering, having $ n $ independent Gumble distributed components with varying location parameter $ \mu. $
\begin{theorem}
Let $ X_{1},....,X_{n}$ $ (Y_{1},....,Y_{n})$  be the set of continuous independent random variables having $ X_{i}\sim Gum(\mu_{i}, \sigma)(Y_{i}\sim Gum(\mu^{*}_{i}, \sigma)),$ $ i=1,2,....,n$. Then, if  $ (\mu_{1},\cdotp\cdotp\cdotp\cdotp\cdotp\mu_{n})\succeq^{m}(\mu^{*}_{1},\cdotp\cdotp\cdotp\cdotp\cdotp\mu^{*}_{n}) $, we have
\begin{center}
	$ X_{1:n}\leq_{hr} Y_{1:n}$.
\end{center}
\begin{proof}
	The  survival function of $ X_{1:n} $ is
	$$\bar{F}_{X_{1:n}}(x)=\prod_{i=1}^{n} \left( 1-e^{-e^{-\frac{(x-\mu_{i})}{\sigma}}}\right) 
	\hspace{.4cm} x,\mu_{i}\in \mathbb{R}, \sigma>0, \hspace{.1cm} i=1,...,n. $$
	It is well known that the sum of the hazard rate functions of the components of a
	series system is equal to the hazard rate function of that system. Therefore,
	we have
	$$ r_{X_{1:n}}=-\frac{d}{dx}\log\left[ \bar{F}_{X_{1:n}}(x)\right]  $$
	
	$$\hspace{1.1cm}=	\frac{1}{\sigma}\sum_{i=1}^{n}\dfrac{e^{-\frac{(x-\mu_{i})}{\sigma}}}{\left[ e^{e^{-\frac{(x-\mu_{i})}{\sigma}}}-1\right] } $$
		$$\hspace{2.5cm}=\frac{1}{\sigma}\sum_{i=1}^{n}\phi(e^{-\frac{(x-\mu_{i})}{\sigma}}), \hspace{.2cm} i=1,...,n.  $$ Where $ \phi(t)=\dfrac{t}{e^{t}-1}. $ Let $ t=e^{-\frac{(x-\mu_{i})}{\sigma}}\in \mathbb{R^{+}},\hspace{.1cm} i=1,...,n.$ Therefore, from Lemma-2.5(1) we say that $ \phi(e^{-\frac{(x-\mu_{i})}{\sigma}}) $ is convex in $ e^{-\frac{(x-\mu_{i})}{\sigma}} $ for $ i=1,...,n. $ Finally, using Lemma-2.6 we conclude that $ r_{X_{1:n}} $ is Schur-convex in $ e^{-\frac{(x-\mu_{i})}{\sigma}} ,$ or equivalently $ r_{X_{1:n}} $ is Schur-convex in $ \mu_{i}$ for $ i=1,...,n. $ Hence, the theorem follows.
	
\end{proof}
\end{theorem}
According to the above result, we immediately obtain the following corollary
\begin{corollary}
		Let $ X_{1},....,X_{n}$ $ (Y_{1},....,Y_{n})$  be the set of continuous  independent random variables having $ X_{i}\sim Gum(\mu_{i}, \sigma)(Y_{i}\sim Gum(\mu^{*}_{i}, \sigma)),$ $ i=1,2,....,n$. Then, if  $ (\mu_{1},\cdotp\cdotp\cdotp\cdotp\cdotp\mu_{n})\succeq^{m}(\mu^{*}_{1},\cdotp\cdotp\cdotp\cdotp\cdotp\mu^{*}_{n}) ,$ and let $ T $ be any random variable which is independent of $ X $ and $ Y $. Then we have
	\begin{center}
		$ X^{T}_{1:n}\leq_{hr} Y^{T}_{1:n}$.
	\end{center}
\begin{proof}
		 By Theorem-2.7(2),  we have to show $ X_{1:n}\leq_{hr} Y_{1:n}$ and $ X_{1:n} $ is DHR for proving $ X^{T}_{1:n}\leq_{hr} Y^{T}_{1:n}$. Now, we know the  hazard rate of $ X_{1:n} $ with respect to $ \mu_{i}, i=1,...,n, $ is 
	$$r_{X_{1:n}}=\frac{1}{\sigma}\sum_{i=1}^{n}\dfrac{e^{-\frac{(x-\mu_{i})}{\sigma}}}{\left[ e^{e^{-\frac{(x-\mu_{i})}{\sigma}}}-1\right] }=\frac{1}{\sigma}\sum_{i=1}^{n}\phi(e^{-\frac{(x-\mu_{i})}{\sigma}}).$$ Where $ \phi(t)=\dfrac{t}{e^{t}-1}. $ Let $ t=e^{-\frac{(x-\mu_{i})}{\sigma}}\in \mathbb{R^{+}},\hspace{.1cm} i=1,...,n.$ Therefore, using Lemma-2.5(2) we say that $ \phi(e^{-\frac{(x-\mu_{i})}{\sigma}}) $ is decreasing in $ e^{-\frac{(x-\mu_{i})}{\sigma}} $ for $ i=1,...,n. $  Equivalently $ r_{X_{1:n}} $ is decreasing in $ \mu_{i}$ for $ i=1,...,n.$ So, we see that $ X_{1:n} $ is DHR.
	Therefore, from Theorem-3.4, and since $ X_{1:n} $ is DHR we conclude that $ X^{T}_{1:n}\leq_{hr} Y^{T}_{1:n},$ this completes the proof.
\end{proof}
\end{corollary}

In the next result, we present the dispersive ordering and the less uncertainty ordering of the series system having independent Gumble distributed components. We find the ordering with respect to the location parameter $ \mu.$
\begin{theorem}
	Let $ X_{1},....,X_{n}$ $ (Y_{1},....,Y_{n})$  be the set of continuous  independent random variables having $ X_{i}\sim Gum(\mu_{i}, \sigma)(Y_{i}\sim Gum(\mu^{*}_{i}, \sigma)),$ $ i=1,2,....,n$. Then, if  $ (\mu_{1},\cdotp\cdotp\cdotp\cdotp\cdotp\mu_{n})\succeq^{m}(\mu^{*}_{1},\cdotp\cdotp\cdotp\cdotp\cdotp\mu^{*}_{n}) $, we have
\begin{enumerate}
	\item $ X_{1:n}\leq_{disp} Y_{1:n}$;
	\item $ X_{1:n}\leq_{LU} Y_{1:n}$.
\end{enumerate}
\begin{proof}
		The the survival function and hazard rate function of $ X_{1:n} ,$ respectively, are
	$$\bar{F}_{X_{1:n}}(x)=\prod_{i=1}^{n} \left( 1-e^{-e^{-\frac{(x-\mu_{i})}{\sigma}}}\right) ,
 $$ and
	$$ r_{X_{1:n}}=\frac{1}{\sigma}\sum_{i=1}^{n}\dfrac{e^{-\frac{(x-\mu_{i})}{\sigma}}}{\left[ e^{e^{-\frac{(x-\mu_{i})}{\sigma}}}-1\right] } 	\hspace{.4cm} x,\mu_{i}\in \mathbb{R}, \sigma>0, \hspace{.1cm} i=1,...,n. $$
	
	  For proving $X_{1:n}\leq_{disp} Y_{1:n},$ and $X_{1:n}\leq_{LU} Y_{1:n},$  by Theorem-2.8 and Theorem-2.9, it is enough if we can show that, under the same condition $ X_{1:n}\leq_{hr} Y_{1:n}$ and $X_{1:n} $ is DHR or 
		$ r_{X_{1:n}} $
		is decreasing in $ \mu_{i}, i=1,...,n,$ hold. Now, from the Theorem-3.3, we observe that  $ X_{1:n}\leq_{hr} Y_{1:n}.$ Next, our claim is $ r_{X_{1:n}} $ is decreasing in $ \mu_{i}, i=1,...,n. $ Now, in the Corollary-3.5, we have already proved that $ r_{X_{1:n}} $ is decreasing with respect to $ \mu_{i}, i=1,...,n, $ i.e. $X_{1:n} $ is DHR, which establish our claim. Thus, the theorem is proved.

\end{proof}

\end{theorem}


\section{Disclosure}
There is no potential conflict of interest and both the authors have equally contributed towards the paper. 
\subsection{Funding}
The research work of Surojit Biswas has been supported by Indian Institute of Technology Kharagpur in the form of M.Tech Teaching Assistantship.

\end{document}